\documentclass[a4paper, reqno, 11pt]{amsart}

\usepackage[english]{babel}
\usepackage{amsmath, amssymb, amsthm, amsfonts, mathtools} 
\usepackage{graphicx, caption, xcolor, placeins} 
\usepackage{enumerate}
\usepackage{ifthen}
\usepackage{bbm}
\usepackage{geometry}
\provideboolean{shownotes} 
\setboolean{shownotes}{true}

\usepackage{color}

\usepackage{tikz} 
\usetikzlibrary{arrows.meta} 
\tikzset{font=\small, 
point/.style={fill, circle, inner sep=1.2pt}, 
>={Straight Barb[round,angle=60:1.2mm 1]} 
} 

\usepackage{hyperref}
\usepackage{setspace}
\usepackage{soul}

\usepackage{xpatch}
\makeatletter
\patchcmd{\@tocline}{\hfil}
{\nobreak\leaders\hbox{\ifnum#1<2\hfill\else$\m@th%
\mkern 4.5 mu\hbox{.}\mkern 4.5 mu$\fi}\hfill\nobreak}{}{}
\def\l@section{\@tocline{1}{10pt}{1pc}{}{\bfseries}}
\def\l@subsection{\@tocline{2}{0pt}{\dimexpr 1pc+2em}{}{}}
\makeatother

\newcommand{\hole}[1]{
\ifthenelse{\boolean{shownotes}}%
{\begin{center} \fbox{ \rule {.25cm}{0cm}
\rule[-.1cm]{0cm}{.4cm} \parbox{.85\textwidth}{\begin{center}
\texttt{#1}\end{center}} \rule {.25cm}{0cm}}\end{center}}
{}
}

\newtheorem{theorem}{Theorem}[section]
\newtheorem{proposition}[theorem]{Proposition}
\newtheorem{lemma}[theorem]{Lemma}

\allowdisplaybreaks
\theoremstyle{remark}

\newcommand{\R}{\mathbb{R}}

\newcommand{\dive}{\mathop{\mathrm {div}}}

\newcommand{\charf}{{{\text{\rm 1}}\kern-.24em {\text{\rm l}}}}

\newcommand{\del}{\partial}
\newcommand{\eps}{\varepsilon}

\newcommand{\cE}{\mathcal{E}}
\newcommand{\supp}{{\rm supp\,}}

\numberwithin{equation}{section}

\begin{document}

\title[Oscillations - Homogenization]{Oscillations in Compressible Navier-Stokes \\ and Homogenization in Phase Transition problems}

\author[A.E. Tzavaras]{Athanasios E. Tzavaras}
\address[Athanasios E. Tzavaras]{
\newline
Computer, Electrical and Mathematical Science and Engineering Division 
\newline
King Abdullah University of Science and Technology (KAUST)
\newline 
Thuwal 23955-6900,  Saudi Arabia
}
\email{athanasios.tzavaras@kaust.edu.sa}

\baselineskip=18pt

\begin{abstract}
In the first part of this article we present some exact solutions for special hyperbolic-parabolic systems with sustained oscillations 
induced by the initial data, most notably the compressible Navier-Stokes system with non-monotone pressure. This part complements 
\cite{Tzavaras23} where such examples are extensively studied. The second part deals with the problem of homogenization for one-dimensional 
models describing phase transitions for viscoelastic materials . Ideas from the kinetic formulation of conservation laws are employed to derive
 effective equations that describe the propagation of oscillations.
\end{abstract}

\maketitle


\section{Introduction}

In this work we study the problem of propagation of oscillations in hyperbolic-parabolic systems.
While certain aspects of this problem are well understood for hyperbolic systems of conservation laws, in one space dimension, 
much less is known at present for hyperbolic-parabolic systems. The issue is already present in the existence theory of compressible Navier-Stokes,
\begin{equation}
\label{intro-compNS}
\begin{aligned}
\del_t \rho + \dive \rho u &= 0
\\
\del_t \rho u + \dive \rho u \otimes u  + \nabla p (\rho) &=  \dive \Big ( \mu (\nabla u + \nabla u^T) + \lambda (\dive u )  \, I \, \Big ) \, ,
\end{aligned}
\end{equation}
where $\rho$ and $u$ are the density and velocity of the fluid while $\mu$ the shear viscosity and $\lambda$ the second viscosity 
satisfy $\mu > 0$, $\lambda + \mu > 0$.
The existence theory of weak solutions of Lions \cite{b-Lions98} and Feireisl \cite{b-Feireisl04} is based on a propagation 
of compactness property from the initial data. Lions conjectured that the compressible Navier-Stokes system may exhibit propagation of
oscillations and gave an example for a forced variant of  \eqref{intro-compNS},  \cite[Rmk 5.8]{b-Lions98}.

Our first goal is to provide an example of propagation of oscillations for the (unforced) compressible Navier-Stokes system \eqref{intro-compNS}
with a non-monotone pressure function.  The example consists of a periodic solution with discontinuous density, defined on a domain $Q = [1,2] \times \R^d$,
\begin{equation}\label{intro-ps}
\begin{aligned}
\rho (t,y) &= 
\begin{cases} 
\; \frac{a}{t^d}  \quad & \qquad \quad k t < |y| < (k + \theta) t \\[5pt]
\; \frac{b}{t^d}  \quad  & \; (k + \theta) t < |y| < (k + 1) t \\
\end{cases}
\; , \qquad k \in \mathbb{N}_0 = \{ 0, 1, 2, ... \} \, ,
\\[5pt]
u(t, y) &= \frac{y}{t} \, .
\end{aligned}
\end{equation}
where $a,b $ be positive constants, $0 < \theta < 1$.
The periodic solution \eqref{intro-ps} produces by space-rescaling an oscillatory solution of compressible Navier-Stokes with sustained oscillations. 
The example generalizes an for the  one-dimensional compressible Navier-Stokes,  \cite{Tzavaras23}, which is, in turn, 
produced as the Eulerian version of a class of oscillating solution examples for viscoelastic models of phase-transitions 
in  Lagrangian coordinates, studied in detail in \cite{KLST23,Tzavaras23}.

The question of calculating effective equations for oscillating solutions of the compressible Navier-Stokes system has been raised and studied 
in \cite{Serre91, Hillairet07}. Serre \cite{Serre91} introduced the problem and provides a detailed study of homogenization 
for the system 
\begin{equation}
\label{intro-vhcg}
\begin{aligned}
w_t - v_x &= 0
\\
v_t - \sigma(w)_x &= ( \frac{\mu}{w} v_x )_x 
\end{aligned}
\end{equation}
for initial data $w(x,0) = w^\eps_0 (x)$, $v(x,0) = v^\eps_0 (x)$. He provided a formal analysis of oscillations for \eqref{intro-compNS}
which was validated rigorously by Hillairet \cite{Hillairet07} under some assumptions on the structure of oscillations. 

We consider here this problem at the level of a viscoelastic model for phase transitions,
\begin{equation}\label{vemodel-intro}
\begin{aligned}
\del_t u &= \del_x v
\\
\del_t v &= \del_x (\sigma (u) + v_x ) \, ,
\end{aligned}
\end{equation}
with non-monotone stress and data $u(x,0) = u_0^\eps(x)$, $v(x,0) = v_0^\eps (x)$.
We introduce to this problem ideas from the kinetic formulation of conservation laws \cite{b-Perthame02}
as it pertains to the study of oscillations for scalar fields \cite{PT00}. A key point is that the weak limits of a uniformly bounded family of functions
$\{ u^\eps \}$ may be represented by the Young measures, but also by the weak-limit of the kinetic function
$$
\charf_{u^\eps < \xi} \rightharpoonup F  \quad \mbox{wk-$\ast$ in $L^\infty$}.
$$
The function $F$ can also be understood as the distribution function of the Young measure $\nu_{t,x} (\xi)$, $\xi \in \R$. We prove in Theorem \ref{mainthm}
that oscillations in solutions $(u^\eps, v^\eps)$ of \eqref{vemodel-intro} are described by the effective system
\begin{equation}\label{intro-effs}
\begin{aligned}
\del_t F &+ \del_\xi \Big ( \big ( v_x + \overline{\sigma(u)} - \sigma(\xi) \big ) F \Big ) + \sigma^\prime (\xi) F = 0   
\\[5pt]
\del_t v &= \del_{xx} v + \del_x  \left (\overline{\sigma(u)} \right )    
\\[5pt]
S &=  \overline{\sigma(u)}  + v_x   = \int \sigma(\xi) dF_{t,x} (\xi)  + v_x  \, .   \\[5pt]
\end{aligned}
\end{equation}
The propagation of $F$ is described by a kinetic equation depending on moments as is typical in the kinetic formulation of systems
of conservation laws, \cite{LPT94, LPS96, PT00}. An alternative (equivalent form) of the effective system appears in section \ref{sec:effs} equations \eqref{effs2}.

The outline of this work is as follows: In section \ref{sec:comprns}, we first review an example from \cite{Tzavaras23} describing
an oscillatory solution for \eqref{intro-vhcg} and the associated solution of its Eulerian counterpart (i.e.,  the system
\eqref{intro-compNS} in one-dimension).  Then we introduce the example \eqref{intro-ps}, prove that it yields a weak solution of
\eqref{intro-compNS} and that it produces via space-rescaling sustainable oscillations. In section \ref{sec:osci}, we study 
oscillatory solutions of the system \eqref{vemodel-intro} (or more precisely the problem  \eqref{vemodel}-\eqref{veid}) induced by oscillations
in the initial data. We prove uniform $L^\infty$ bounds for the sequence $\{ u^\eps \}$ using a variant of an  argument from \cite{AB82};
then prove compactness for the sequence $\{S^\eps\}$ using parabolic theory. 
After recalling the relevant aspects of the theory of kinetic functions \cite{b-Perthame02},
they are used in order to calculate the effective equations.

%
%
%


\section{Sustained oscillations for compressible Navier-Stokes}\label{sec:comprns}

In this section we outline examples of oscillating solutions for the compressible Navier-Stokes system 
\eqref{intro-compNS}. A first example  \cite[sec 7]{KLST23} concerns the system \eqref{vemodel-intro} describing shear motions
for viscoelasticity of the rate-type with non-monotone stress. This example will not be presented in detail here.
However, in section \ref{sec:osci},  we study the homogenization properties of \eqref{vemodel-intro}.

\subsection{Sustained oscillations for viscous gases in one dimension}\label{sec:1dnonmon}

By contrast, we will present an outline of a closely related example from \cite{Tzavaras23} for the system
\begin{equation}
\label{longonev}
\begin{aligned}
w_t &= v_x
\\
v_t &= \sigma(w)_x +  \del_x \left ( \frac{\mu}{w} v_{x} \right )
\end{aligned}
\end{equation}
describing longitudinal motions  $y(t,x) :(0,T) \times [0,1] \to \R$ for a bar, with $u = y_x > 0$
the longitudinal strain and $v= y_t$ the velocity. The total stress 
$$
S = \sigma (y_x) +    \frac{\mu }{y_x} y_{t x}
$$ 
has an elastic component $\sigma_{el} = \sigma (w)$ and a viscous component $\sigma_{v} = \frac{\mu}{w} v_x$ with viscosity
$\frac{\mu}{w}$. The system \eqref{longonev} can also be viewed as a description of isothermal (or barotropic) gas-dynamics in Lagrangian coordinates
for a viscous gas; then $\sigma (w) = - p(w)$ with $p$ the gas pressure.

Assume two positive states $0 < a < b$ are fixed such that $0 < a < 2a < b < 2b$ 
and suppose the stress function $\sigma (w)$ satisfies
\begin{equation}
\label{condnonm}
 \sigma (\tau  a) = \sigma ( \tau  b) \qquad \mbox{ for $\tau \in [1,2]$} \, .
\end{equation}
The function $\sigma(w)$ is necessarily non-monotone.
The motion $Y(t,x)$ is defined for $(t,x) \in [1,2]\times \R$, as follows: 
The function $W$ is selected to be periodic
\begin{equation}\label{perU}
W(t,x)  := \begin{cases}
 a \,  t   &  \quad  \; \;  k <  x <  k + \theta
 \\
 b \,  t   &  k + \theta < x < k + 1
 \end{cases}
\quad \; k \in \mathbb{Z}  \, ,
\end{equation}
and, for $c_\theta = \theta a + (1-\theta) b$,  define
\begin{equation}
\label{perY}
\begin{aligned}
Y(t,x) = \int_0^x W (t, y) dy &= 
\begin{cases}
k c_\theta t + (x-k) a t  & \quad  \; \;  k <  x <  k + \theta
\\
 k c_\theta t + \theta a t + \big ( x - k - \theta) b t   &  k + \theta < x < k + 1
\end{cases} \, , \quad k \in \mathbb{Z} \, ,
\\
V(t,x) = \del_t Y (t,x) &=
\begin{cases}
k c_\theta  + (x-k) a   & \quad  \; \;  k <  x <  k + \theta
\\
 k c_\theta  + \theta a  + \big ( x - k - \theta) b    &  k + \theta < x < k + 1
\end{cases} \, , \quad k \in \mathbb{Z} \, .
\end{aligned}
\end{equation}
Then $(W, V)$ is a weak solution of \eqref{longonev} on $(1,2) \times\R$, the equations are satisfied 
in a classical sense on $(1,2)\times (k, k + \theta)$, $(1,2)\times(k + \theta , k+1)$. 
Following \cite{Hoff86}, the Rankine-Hugoniot conditions at the steady interfaces $x = k$ and  $x = k + \theta$, $1 \le t \le 2$,
reduce to the condition
\begin{equation}\label{sec7RH}
s= 0 \, , \quad \left [ \sigma (w) +  \mu \frac{v_x }{w} \right ] = 0 \, ,
\end{equation}
which is fulfilled due to \eqref{condnonm}.  We refer to \cite{KLST23, Tzavaras23} for details of the construction.

The function $Y(t,x)$ is  rescaled and restricted to $(t,x) \in  Q = (1,2)\times (-1,1)$ by
$$
y_n (t,x) = \frac{1}{n} Y(t, nx)  \, \quad w_n = \del_x y_n = W  (t, nx)  \, \quad {v_n} = \del_t y_n  =  \frac{1}{n}  V (t, nx)
$$
Then $(w_n , v_n )$ is a weak solution of \eqref{longonev} and we compute the limits
$$
\begin{aligned}
w_n  &\rightharpoonup  (a\theta + b(1-\theta) ) t  \quad \mbox{ weakly-$\star$ in  $L^\infty \big ( Q \big )$ }
\\
\del_x v_n &\rightharpoonup  (a\theta + b(1-\theta) )  \quad \mbox{ weakly-$\star$ in  $L^\infty \big ( Q \big )$ }
\\
v_n &\to (a\theta + b(1-\theta)) x  \quad \mbox{ strongly in $L^2 \big ( Q \big )$ }
\\
\sigma (w_n) &\rightharpoonup  \theta \sigma ( at) + (1-\theta) \sigma(bt) \quad \mbox{ weak-$\star$ in $L^\infty(Q)$ }
\end{aligned}
$$
The oscillations of the sequence $\{ w_n \}$  are described  by the associated Young measure 
$$\nu = \theta \delta_{at} + (1-\theta) \delta_{b t}$$
and are induced by oscillations in the initial data $w_n (1,x)$ and $\del_x v_n (1, x)$.

Next, consider the one-dimensional version of the compressible Navier-Stokes system,
\begin{equation}
\label{compNS1d}
\begin{aligned}
\rho_t  + (\rho u)_y  &= 0
\\
(\rho u)_t  + ( \rho u^2  +  p (\rho) )_y  &=  \mu u_{y y}
\end{aligned}
\end{equation}
Here, $\rho (t,y)$ and $u(t,y)$ are the density and velocity of the gas in Eulerian coordinates $(t,y)$ while $\mu > 0 $ is the viscosity.
The system \eqref{compNS1d}  is the Eulerian version of   \eqref{longonev}, obtained through the transformation
\begin{equation}\label{transf}
\begin{aligned}
\frac{\del y}{\del t} (t,x) &= v (t, x) = u (t, y(t,x)) \, ,
\\
 \rho (t, y(t,x))  &= \frac{1}{\frac{\del y}{\del x} (t,x)} = \frac{1}{w(t,x)}
\end{aligned}
\end{equation}
where $p(\rho) =  - \sigma (w)$.

When  \eqref{perU} and \eqref{perY} are transformed to Eulerian coordinates, via \eqref{transf},
they provide a weak solution for \eqref{compNS1d} on the interval $(1.2) \times \R$:
\begin{align}
u(t,y) &= \frac{y}{t}
\label{eqnu}
\\
\rho(t,y) &= 
\begin{cases} \;  \frac{1}{ta}  \quad & \quad  0 < \frac{y}{t}  - k v_0 (\theta) <  a \theta   \\[5pt]
                      \;  \frac{1}{tb}  \quad &  \; \; a \theta < \frac{y}{t}  - k v_0 (\theta) < v_0 (\theta)
\end{cases}
\quad k \in \mathbb{Z}
\label{eqnrho}
\end{align}
see \cite[sec 5]{Tzavaras23}. This solution is valid on the domain $(1.2) \times \R$ under the restriction for the pressure
\begin{equation}\label{hyppr1d}
p \left ( \frac{1}{at} \right ) = p \left ( \frac{1}{bt} \right )  \quad \mbox{for \;  $t \in [1,2]$}.
\tag {{H}$_{1d}$}
\end{equation}
By rescaling it, 
\begin{equation}\label{oscisol}
u_n (t, y) = \frac{1}{n} u (t, n y) \, , \quad \rho_n (t, y) = \rho (t, n y)   \qquad y \in (-1,1), \;  t \in [1,2] \, ,
\end{equation}
we obtain a sequence  of weak solutions $(\rho_n , u_n )$ that exhibit sustained oscillations and whose  weak limit $(\bar \rho, \bar u)$ and $\bar p$
does not (in general) satisfy the relation $p ( \bar \rho) = \bar p$.

\subsection{A weak solution for compressible Navier Stokes with non-monotone pressure}
Consider first the pressureless Euler equations in $d$ dimensions,
\begin{equation}
\label{pleuler}
\begin{aligned}
\del_t \rho + \dive \rho u &= 0
\\
\del_t \rho u + \dive \rho u \otimes u   &=  0
\end{aligned}
\end{equation}
A direct computation shows that for $\rho_0 > 0$ constant, the function
\begin{equation}\label{unif}
\hat \rho = \frac{\rho_0}{t^d} \, , \quad \hat u = \frac{y}{t}  \, , \qquad y \in \mathbb{R}^d \, , \; t > 0 \, ,
\end{equation}
is an exact solution of \eqref{pleuler}. Indeed,
$$
\begin{aligned}
\left ( \frac{\del}{\del t} + \hat u_j \frac{\del}{\del y_j} \right ) \frac{y_i}{t} &= 0
\\[5pt]
\del_t \frac{\rho_0}{t^d}   + \dive \left ( \frac{\rho_0}{t^d} \frac{y}{t}  \right ) & = 0
\end{aligned}
$$

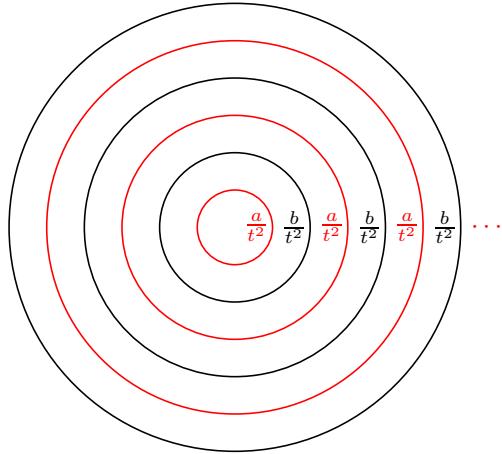
\begin{figure}
\begin{tikzpicture}[scale=1.1] 
\def\nn{6} 
\foreach \i in {1,3,...,\nn} { 
	\draw[red, semithick] (.45*\i,0) arc (0:360:.45*\i) node[xshift=-6]{$\frac{a}{t^2}$}; 
	\draw[semithick] ({.45*(\i+1)},0) arc (0:360:{.45*(\i+1)}) node[xshift=-6]{$\frac{b}{t^2}$}; 
}
\node[red] at (.45*\nn,0)[right] {$\cdots$}; 
\end{tikzpicture} 
\caption{ density jumps at moving interfaces - dimension $d=2$}\label{fig:1}
\end{figure}

Let now $a,b $ be positive constants, $0 < \theta < 1$ and consider the function $(\rho, u)$ defined on the
domain $Q = [1,2] \times \R^d$ by
\begin{equation}\label{propsol}
\begin{aligned}
\rho (t,y) &= 
\begin{cases} 
\frac{a}{t^d}  \quad & \qquad \quad k t < |y| < (k + \theta) t \\[5pt]
\frac{b}{t^d}  \quad  & \; (k + \theta) t < |y| < (k + 1) t \\
\end{cases}
\; , \qquad k \in \mathbb{N}_0 = \{ 0, 1, 2, ... \} \, ,
\\[5pt]
u(t, y) &= \frac{y}{t} \, .
\end{aligned}
\end{equation}
In \eqref{propsol}, $u$  is smooth and $\rho$ has discontinuities  at interfaces depicted in Figure \ref{fig:1}.
Note that $(\rho, u)$ is an exact solution in each of the annuli $k t < |y| < (k + \theta) t$ and $(k + \theta) t < |y| < (k + 1) t$, $k \in \mathbb{N}_0$.

The discontinuities move on the shock surfaces $S_k$ defined by $\xi_k (t) = k t \hat y$, $\hat y  = \frac{y}{|y|} \in \mathcal{S}^{d-1}$ 
and on the shock surfaces $S_{k + \theta}$ defined by $\xi_{k + \theta}  (t) = (k+\theta)  t \hat y$, $\hat y  = \frac{y}{|y|} \in \mathcal{S}^{d-1}$.
On the shock surface $S_k$ the Rankine-Hygoniot conditions are satisfied:
$$
\begin{aligned}
-s [\rho] + \nu_j [\rho u_j ] &= (\nu_j u_j - s ) [\rho] = (\hat y \cdot u |_{S_k} -  |\dot \xi_k | )   [\rho] = 0
\\
-s [\rho u_i ] + \nu_j [\rho u_i u_j ] &= u_i  (\hat y \cdot u |_{S_k} -  |\dot \xi_k | )  [\rho] = 0
\end{aligned}
$$
where the summation convention over repeated indices was used.
A similar calculation shows that the Rankine-Hugoniot conditions are also satisfied on $S_{k + \theta}$. We conclude that $(\rho, u)$
defines a weak solution of \eqref{pleuler} on the domain $\{ t >0 \} \times \R^d$.

Consider next the compressible Navier-Stokes system
\begin{equation}
\label{compns}
\begin{aligned}
\del_t \rho + \dive \rho u &= 0
\\
\del_t \rho u + \dive \rho u \otimes u  + \nabla p (\rho) &=  \dive \Big ( \mu (\nabla u + \nabla u^T) + \lambda (\dive u )  \, I \, \Big )
\end{aligned}
\end{equation}
where $\lambda$, $\mu$ are constants with $\lambda + \mu > 0$, $\mu > 0$. Observe that $(\hat \rho, \hat u)$ in \eqref{unif} is still
an exact solution of \eqref{compns} for any pressure function $p (\rho)$. Consider now the function $(\rho, u)$ defined by \eqref{propsol}
and test under what conditions it is a weak solution on the domain $(1,2) \times \R^d$. We examine the Rankine-Hugoniot conditions on
the surfaces $S_k$ and $S_{k+\theta}$. Note that the balance of mass is automatically satisfied and the balance of momentum reduces to
$$
n_j  \left [ p(\rho) \delta_{i j} +  \mu \left(  \frac{\del u_i}{\del y_j}  + \frac{\del u_j}{\del y_i} \right )  + \lambda \frac{\del u_k }{\del y_k} \delta_{i j}  \right ] = 0
$$
In turn, this reduces to satisfying the condition
\begin{equation}\label{hypprmd}
p \left ( \frac{a}{t^d} \right ) = p \left ( \frac{b}{t^d} \right )  \quad \mbox{for \;  $t \in [1,2]$}.
\tag {{H}$_{md}$}
\end{equation}
This can be achieved by selecting $0 < a<b$ to satisfy the strengthened condition $a < \frac{b}{2^d}$ and then selecting $p(\rho)$ so that \eqref{hypprmd}
is satisfied. We summarize:

\begin{proposition}
Let the pressure satisfy \eqref{hypprmd} for some $0 < \frac{a}{2^d} < a < \frac{b}{2^d} < b$. Then, \eqref{propsol} is a weak solution,
defined for $y \in \R$,  $t \in [1,2]$, 
of the compressible Navier-Stokes system \eqref{compns} with non-monotone pressure.  
\end{proposition}

\subsection{Sustained oscillations}
The solution \eqref{propsol} will be rescaled  to provide sustained oscillations as follows:
\begin{equation}\label{oscisol2}
\begin{aligned}
u_n (t, y) &= \frac{1}{n} u (t, n y) = \frac{y}{t}  \\
 \rho_n (t, y) &= \rho (t, n y)   
\end{aligned}
\quad , \qquad   t \in [1,2] \, , \quad y \in \R^d \\[5pt]
\end{equation}
where $(\rho, u)$ is given by \eqref{propsol}, $n \in \mathbb{N}$. 

Note the direct calculations
\begin{align*}
\del_t \rho_n + \dive \rho_n u_n &=  \big ( \del_t \rho + \dive \rho u \big ) \Big |_{(t, ny)}
\\
\del_t  (\rho_n u_n)  + \dive ( \rho_n u_n \otimes u_n ) &=  \frac{1}{n} \big ( \del_t  (\rho u)  + \dive ( \rho u \otimes u \big ) \Big |_{(t, ny)}
\\
 \nabla p(\rho_n) + \dive \Big ( \mu (\nabla u_n + \nabla u_n^T) + \lambda (\dive u_n )  \, I \, \Big ) &= n  \left [  \nabla p(\rho) 
 + \dive \Big ( \mu (\nabla u + \nabla u^T) + \lambda (\dive u )  \, I \, \Big ) \right ] \Big |_{(t, ny)}
 \end{align*}
 Therefore, although the equation \eqref{compns} is not invariant under the scaling \eqref{oscisol}, because $(\rho, u)$ in \eqref{propsol} satisfies
 simultaneously the pressureless Euler equation \eqref{pleuler}, and the steady equation
 $$
  \nabla p(\rho) + \dive \Big ( \mu (\nabla u + \nabla u^T) + \lambda (\dive u )  \, I \, \Big ) = 0
 $$
 it follows that $(\rho_n, u_n)$ also satisfies the compressible Navier Stokes system.
 
 When the function $(\rho_n, u_n)$ is restricted to the domain $Q = [1,2] \times B_1$ it provides $u_n = \frac{y}{t}$ and
 \begin{equation}\label{propsol2}
\begin{aligned}
\rho_n (t,y) &= 
\begin{cases} 
\frac{a}{t^d}  \quad &  \qquad \frac{k}{n}   < \frac{|y|}{t}  < \frac{k + \theta}{n}   \\[5pt]
\frac{b}{t^d}  \quad & \quad  \frac{k + \theta}{n} < \frac{|y|}{t}  < \frac{k+1}{n}  \\
\end{cases}
\; , \qquad   k = 0, 1, ... , n-1  \, .
\end{aligned}
\end{equation}
$(\rho_n, u_n)$ are weak solutions of \eqref{compns} with the properties $u_n \to u = \frac{y}{t}$ strongly while
$$
\rho_n \rightharpoonup \rho = \theta \frac{a}{t^d} + (1-\theta) \frac{b}{t^d} \quad \mbox{ weakly-$\star$ in  $L^\infty \big ( Q \big )$ }
$$
 The Young measure associated with the sequence $\rho_n$ is 
 $$
 \nu = \theta \delta_{\frac{a}{t^d}} + (1-\theta) \delta_{\frac{b}{t^d}}
 $$
 The sequence $\{ \rho_n\}$ exhibits persistent oscillations tinduced  by oscillations in the initial data.
 While $(\rho, u)$ is a solution of compressible Navier Stokes, the reader should note that 
 $\{ \rho_n \}$ has the property that $p(\rho_n)  \rightharpoonup q$ with  in general $q  \ne p(\rho)$.


\section{Homogenization for one-dimensional viscoelastic models with phase transitions.} \label{sec:osci}

We consider next the problem of calculating effective equations for the propagation of oscillations for models
describing phase transitions for rate-dependent materials in one space dimension. The system
\begin{equation}\label{vemodel}
\begin{aligned}
\del_t u &= \del_x v
\\
\del_t v &= \del_x (\sigma (u) + v_x ) \, ,
\end{aligned}
\end{equation}
describes shear motions of a viscoelastic material, with $u$ the shear deformation, $v$ the velocity in the shear direction,
and $S$ the total stress
\begin{equation}\label{defstress}
S = \sigma (u) + v_x
\end{equation} 
all defined for $0 \le x \le 1$, $t > 0$. We take up the simplest case of 
traction-free boundary conditions
\begin{equation}\label{vebc}
S (0,t) = S(1,t) = 0
\end{equation}
with highly oscillatory initial data
\begin{equation}\label{veid}
u^\eps (x, 0)= u_0^\eps (x) \, , \quad v^\eps (x, 0) = v_0^\eps (x) \, ,
\end{equation}
that are assumed to satisfy for some $p>1$  the uniform bounds
\begin{equation}\label{ubdata}
\| u^\eps_0 \|_{L^p} \le C \, , \quad  \| v^\eps_0 \|_{L^2} \le C   \tag{A$_1$}.
\end{equation}
Oscillatory data in general induce oscillations to solutions
$(u^\eps, v^\eps)$ of \eqref{vemodel}-\eqref{veid}. Our scope is to compute the effective equations describing the oscillations.
In the sequel, we drop the $\eps$-dependence where possible, but it will always be implied.

The stress function is non-monotone as appropriate for models of phase transitions with behavior $\sigma(u) \to \pm \infty$ as $u \to \pm \infty$. The stored energy
$W(u) = \int_0^u \sigma (\tau) d\tau$ is assumed to obey for some $p > 1$ the bounds
\begin{equation}\label{boundW}
c | u|^p - C \le W(u) \le C |u|^p   \tag{H$_1$}
\end{equation}
for some $c, C  > 0$ constants.

The existence theory of the viscoelasticity model \eqref{vemodel} is well understood and the reader is referred to \cite{AB82}, \cite{KLST23} and references therein. 
For initial data uniformly bounded,  $u_0 \in L^p$, $p > 1$ and  $v_0 \in L^2$, the formal energy identity
$$
\del_t \Big ( \frac{1}{2} v^2 + W(u) \Big ) - \del_x \big (v (\sigma (u) + v_x ) \big ) + v_x^2 = 0
$$
yields the uniform bounds
\begin{equation}\label{energyb}
\int \tfrac{1}{2} (v^\eps)^2 (x,t)  + W(u^\eps (x,t) ) \, dx   + \int_0^t \int (v^\eps_x)^2 \le \int \tfrac{1}{2} (v_0^\eps)^2 + W( u_0^\eps ) \, dx = \cE^\eps (0)  \le C
\end{equation}
Using \eqref{boundW}, the standard theory of weak convergence, and the Aubin-Lions Lemma, we conclude that in the framework of the energy bounds \eqref{energyb}, we have the convergence information:
\begin{itemize}
\item{(i)} $u^\eps \rightharpoonup u  \quad \mbox{wk-$\ast$  in $L^\infty (L^p)$}$ 
\item{(ii)} $v^\eps \to v  \quad  \mbox{strongly in  $L^2 (L^2)$}$ 
\end{itemize}
along subsequences if necessary.

\subsection{$L^\infty$ bounds for the strain}
Next, we improve the a-priori estimates under additional uniform bounds on the initial data:
\begin{equation}\label{ub2data}
\| u^\eps_0 \|_{L^\infty} \le C \, , \quad  \| v^\eps_0 \|_{W^{1,\infty}} \le C   \tag{A$_2$}.
\end{equation}
An idea of Andrews-Ball \cite{AB82} can be adapted  to yield uniform bounds on the shear strain. This
hinges on a technical hypothesis on the stress function (see \eqref{hypgr} below).

Due to the stress-free boundary condition \eqref{vebc}, we obtain from \eqref{vemodel}, \eqref{energyb},
$$
\del_t \int_0^x v(y,t) dy = S(x,t) - S(0,t) = \sigma (u(x,t))  + u_t (x,t) \, ,
$$
with
\begin{equation}\label{momb}
\Big | \int_0^x v(y, t) dy \Big | \le \Big ( \int_0^1 v^2(y,t) dy \Big )^{\tfrac{1}{2}} \le \sqrt{\cE(0)} = : M
\end{equation}
This yields the identity
$$
\del_t \Big ( u(x,t) - \int_0^x v(y,t) dy \Big ) + \sigma (u(x,t)) = 0
$$
which after some rearrangements is rewritten as
$$
\begin{aligned}
\del_t \Big ( u(x,t) - \int_0^x v(y,t) dy \Big )  &+ \sigma \Big ( u(x,t) - \int_0^x v(y,t) dy \Big ) 
\\
&= \int_0^1 \frac{d}{ds} \sigma \Big ( u(x,t) -  s \int_0^x v(y,t) dy \Big ) \, ds
\\
&= - \left ( \int_0^1 \sigma^\prime \Big ( u(x,t) -  s \int_0^x v(y,t) dy \Big ) \, ds \right ) \int_0^x v(y,t) dy \, .
\end{aligned}
$$

In summary, setting
\begin{equation}\label{defg}
g(x,t) = u(x,t) - \int_0^x v(y,t) dy  \, , \quad  A(x,t) = \int_0^x v(y,t) dy \, ,
\end{equation}
we conclude that $g$ satisfies the ordinary differential equation
\begin{equation}\label{eqng}
\del_t g + \sigma (g) = - \left (\int_0^1  \sigma^\prime \big ( g(x,t) + (1-s) A(x,t) \big ) ds \right ) A(x,t)
\end{equation}
Moreover, by virtue of \eqref{momb}, $g$ also satisfies the differential inequalities
\begin{align}
\label{ineq1}
\del_t g + \sigma (g)  &\le \max_{s \in [0,1] \, \; |A| \le M} \big | \sigma^\prime ( g + (1-s) A ) \big | | A|  = : B(g, M)
\\
\del_t g + \sigma (g)  &\ge  - \max_{s \in [0,1] \, \; |A| \le M} \big | \sigma^\prime ( g + (1-s) A ) \big | | A|  =   - B(g, M)
\end{align}

\begin{lemma} \label{lem1}
Let 
$$
B(u,M) = \max_{s \in [0,1] \, \; |A| \le M} \big | \sigma^\prime ( u + (1-s) A ) \big | | A| 
$$ 
and suppose that $\sigma (u)$ has the limiting behavior
\begin{equation}\label{hypgr}
\begin{aligned}
\lim_{u \to \infty} \big ( \sigma (u) - B(u,M) \big ) &= + \infty
\\
\lim_{u \to - \infty} \big ( \sigma (u) + B(u,M) \big ) &= - \infty
\end{aligned}
\tag {H$_2$}
\end{equation}
If the initial data satisfy the bounds \eqref{ubdata}, \eqref{ub2data} then for solutions $(u^\eps, v^\eps)$ of \eqref{vemodel}-\eqref{veid}
\begin{equation}\label{ubound}
| u^\eps (x,t) | \le K
\end{equation}
where $K$ is a constant independent of $\eps$.
\end{lemma}

\noindent
Observe that, when $\sigma (u)$ grows like a power $u^{p-1}$ as $u \to \infty$ with $p>2$ then $\sigma^\prime (u)$ grows like $u^{p-2}$ and 
thus \eqref{hypgr} is satisfied as $u \to \infty$; similarly for $u \to -\infty$.

{\it Proof.} 
Suppose $g$ satisfies \eqref{eqng} and let $|g_0 (x) | \le K$.  By virtue of \eqref{hypgr} we may select $g_+ > 0$, $g_- < 0$ such that
$$
\begin{aligned}
\sigma(g ) - B(g, M) > 0  \quad \mbox{for $g > g_+$}
\\
\sigma(g ) + B(g, M) < 0  \quad \mbox{for $g < g_-$}
\end{aligned}
$$

Then we claim that
$$
\min\{ -K, g_- \} < g (x,t) < \max\{K, g_+ \}  \, .
$$
Indeed, by the assumption of the data, $- K < g(x,t) < K$ for $t$ sufficiently small. Let $t$ be the first time that for some $x$ we have $g(x,t) = \max\{K, g_+ \} $.
Then \eqref{ineq1} implies that at $(x,t)$ as above it is $g_t (x,t) \le 0$ which gives a contradiction. Similarly, follows the lower bound.

Since $A(x,t)$ satisfies \eqref{momb}, equation \eqref{defg} implies that $|u(x,t)| \le K$. The bound under hypotheses \eqref{ubdata},
\eqref{ub2data} is uniform in $\eps$.
\qed

\subsection{Regularizing effect for the total stress}
From \eqref{vemodel} follows that the total stress  $S = \sigma (u) + v_x$
satisfies the linear parabolic equation
\begin{equation}\label{eqstress}
\del_t S = \del_{x x} S + \sigma^\prime (u) \big ( S - \sigma (u) \big )
\end{equation}
with boundary conditions \eqref{vebc} and initial conditions
$$
S (x,0) = \sigma (u_0^\eps) + v_{0 x}^\eps \in L^\infty \, ,
$$
uniformly bounded by  \eqref{ub2data}.

\begin{lemma}\label{regstress}
Under the hypotheses of Lemma \ref{lem1} the total stress $S^\eps$ satisfies the uniform bound
\begin{equation}\label{stressb}
| S^\eps (x,t) | \le C 
\end{equation}
for $(x,t) \in Q_T = (0,1) \times (0,T]$, and along a subsequence
\begin{equation}\label{stressconv}
S^\eps \to S \quad \mbox{ a.e. $(x,t) \in Q_T$}.
\end{equation}
\end{lemma}

{\it Proof.} The stress $S^\eps$ satisfies a boundary value problem
\begin{equation}
\begin{cases}
\del_t S = \del_{xx} S + a^\eps (x,t) S + b^\eps (x,t)  \quad (x,t) \in Q_T  & \\[3pt]
S(0,t) = S(1,t) = 0   & \\[3pt]
S(x,0) = S_0^\eps (x)  & \\
\end{cases}
\end{equation}
with coefficients $a^\eps$, $b^\eps$ uniformly bounded in $L^\infty(Q_T)$. Standard parabolic theory for linear equations
\cite[Thm 7.1]{b-LSU68}  implies the $L^\infty$ uniform bound. In addition, Theorem \cite[Thm 10.1]{b-LSU68} implies that
there is an exponent $0 < \alpha < 1$ so that $S(x,t)$ is H\"{o}lder continuous of exponents $\alpha$-$\alpha/2$. More precisely, 
for any parabolic cylinder $Q^\prime \subset Q_T$ separated by a distance $d > 0$ from
the parabolic boundary $\Gamma_T$ the H\"{o}lder-seminorm 
$$
< S^\eps >_{\alpha , \alpha/2} \le M_d
$$
with the constant depending on $d$ and blowing-up as $d$ tends to zero. Since the boundary data are H\"{o}lder continuous,
the estimate holds for parabolic cylinders that reach the lateral side of $\Gamma_T$ but stay at a distance $d>0$ from the initial plane $x=0$ (see 
second part of \cite[Thm 10.1]{b-LSU68}).
On each cylinder $Q^\prime$, the Ascoli-Arzela implies that there is $S \in C^{\alpha, \alpha/2} (Q^\prime)$ such that along a subsequence $S^\eps \to S$ pointwise on $Q^\prime$.
Taking a sequence of approximations with $d_k$ the distance to the initial plane tending to zero and using a diagonal
argument,  we conclude that there is a subsequence $S^\eps$ and $S$ so that \eqref{stressconv} holds.  In fact, the convergence $S^\eps \to S$ 
is that of uniform convergence on compact subsets of $[0,1] \times (0, T]$. \qed

\subsection{The effective system - Preliminaries}
Consider next a family $(u^\eps, v^\eps)$ of solutions to \eqref{vemodel}-\eqref{veid} and proceed to calculate the effective
equation describing the homogenization process. Under hypotheses \eqref{ub2data}  the family $(u^\eps, v^\eps)$ satisfies \eqref{ubound}, \eqref{stressb} and has the convergence properties
that, along a subsequence,
\begin{equation}\label{wklim}
\begin{aligned}
u^\eps \rightharpoonup u  \quad \mbox{wk-$\ast$  in $L^\infty (Q_T)$}
\\
v^\eps \to v  \qquad \qquad  \mbox{in  $L^2 (Q_T)$}
\\
S^\eps \to S \quad \mbox{ a.e. $(x,t) \in Q_T$}
\end{aligned}
\end{equation}
We introduce the Young measure associated with the sequence $\{ u^\eps\}$, namely a parametrized family of probability measures
$\nu_{t,x} (\lambda)$ that represents weak limits
\begin{equation}\label{defYM}
\beta(u^\eps) \rightharpoonup \langle \nu_{t,x} , \beta(\lambda) \rangle  = \int \beta (\lambda) \, d\nu_{t,x} (\lambda)  \qquad \mbox{wk-$\ast$  in $L^\infty (Q_T)$}
\end{equation}
for any $\beta(u)$ continuous. Due to the uniform bound \eqref{ubound} the support of the Young measure $\supp \nu_{t,x} \subset [-K, K]$.

We start along the lines of \cite{Serre91} and compute 
\begin{equation}\label{epsfamily}
\del_t \beta (u^\eps) = (S^\eps - \sigma(u^\eps) ) \beta^\prime (u^\eps)
\end{equation}
and using the Young measure representation and Lemma \ref{regstress} we derive the family of equations
\begin{equation}\label{family}
\del_t \Big\langle \nu_{t,x} , \beta(\lambda) \Big\rangle = \Big  \langle \nu_{t,x} ,  S \beta^\prime (\lambda)  - \sigma (\lambda) \beta^\prime (\lambda)  \Big \rangle 
\end{equation}
The reader is referred to  \cite{Serre91} for the presentation of an effective equation representing the family of equations \eqref{family}.
The effective equation is formal and is based  on a rendition of the Young-measures (for scalar functions) via 
rearrangements, an idea attributed to  L. Tartar.

Here, we propose an alternative way of representing the family \eqref{family} using the  kinetic function as a tool to describe propagations of oscillations
(see \cite[Ch 2]{b-Perthame02} and \cite{PT00}). The bound \eqref{ubound} implies for $\xi \in \R$
\begin{equation}\label{defF}
\charf_{u^\eps < \xi}  \rightharpoonup F(t,x,\xi) \qquad \mbox{wk-$\ast$  in $L^\infty (Q_T)$}
\end{equation}
to some $F \in L^\infty (Q_T \times \R)$.  The function $F$ is monotone increasing in $\xi$ and satisfies $F(t,x, -\infty) = 0$, $F(t,x,\infty)=1$.
 The reader will understand the role of the function $F$ by noting that for $\beta \in C^1_c (\R)$
we have
\begin{equation}\label{repform1}
\begin{aligned}
\beta(u^\eps) &= - \int_{u^\eps}^{\infty} \beta^\prime (\xi) d\xi  
\\
&= - \int_{-\infty}^{\infty} \beta^\prime (\xi) \charf_{u^\eps < \xi} d\xi  \rightharpoonup - \int_{-\infty}^{\infty} \beta^\prime (\xi) F(t,x,\xi) d\xi
\end{aligned}
\end{equation}
\begin{equation}\label{repform2}
\begin{aligned}
\beta(u^\eps) &=  \int^{u^\eps}_{-\infty} \beta^\prime (\xi) d\xi  
\\
&= \int_{-\infty}^{\infty} \beta^\prime (\xi) \charf_{\xi < u^\eps} d\xi  \rightharpoonup  \int_{-\infty}^{\infty} \beta^\prime (\xi) (1 - F(t,x,\xi)) d\xi
\end{aligned}
\end{equation}
Comparing these relations with \eqref{defYM} we see that
\begin{equation}\label{ibp}
\int \beta (\lambda) \, d\nu_{t,x} (\lambda) = - \int_{-\infty}^{\infty} \beta^\prime (\xi) F(t,x,\xi) d\xi =  \int_{-\infty}^{\infty} \beta^\prime (\xi) (1 - F(t,x,\xi)) d\xi
\quad \mbox{for $\beta \in C^1_c(\R)$}.
\end{equation}
Hence, in the sense of distributions
$$
\nu_{t,x} (\xi) = \del_{\xi } F(t,x,\xi) \, .
$$
An alternate way to introduce $F(t,x,\xi)$ is to realize that since $\nu_{t,x} (\lambda)$ is a Borel measure on $\R$ one may introduce its distribution function
\begin{equation}\label{defdistr}
F(t,x,\xi) = \int_{(-\infty, \xi]} d\nu_{t,x}(\lambda)
\end{equation}
This choice picks the right continuous representative of the function $F_{t,x} (\xi)$. Then  the integration by parts formula 
\cite[Thm 3.30]{b-Folland99} for BV functions provides an alternate viewpoint to \eqref{ibp}.

We will use the notation $dF_{t,x} (\xi)$ for integration against the Young measure. Then the weak$-\ast$ limit in $L^\infty$ of a continuous function $f(u)$
can be computed via the formula
\begin{equation}\label{moments}
\begin{aligned}
\overline{f(u)} &= \int f(\lambda) d\nu_{t,x} (\lambda) = \int f(\xi) dF_{t,x} (\xi)
\\
&= \int_{(-\infty, 0]} f(\xi) dF_{t,x} (\xi) + \int_{(0, \infty)} f(\xi) d(F_{t,x} (\xi) - 1)
\\
&= f(0) - \int_{-\infty}^0 f^\prime (\xi) F(t,x,\xi) d\xi - \int_0^\infty f^\prime (\xi) ( F(t,x,\xi) - 1) d\xi 
\end{aligned}
\end{equation}
where we have used the integration by parts formula, with the right continuous representative of $F(\xi)$, and the facts that $\nu$ has compact support 
and $F(-\infty) =0$, $F(+\infty) =1$.
Formula \eqref{moments} is used for computing the weak limits and should be compared to corresponding formulas using the kinetic function in 
\cite[Ch 2]{b-Perthame02}.

\subsection{Computation of the effective system}\label{sec:effs}
We compute two forms of effective systems describing the propagation of oscillations in \eqref{vemodel}-\eqref{veid}.
The first system has the form
\begin{equation}\label{effs1}
\begin{cases}
\del_t F + \del_\xi \Big ( \big ( v_x + \overline{\sigma(u)} - \sigma(\xi) \big ) F \Big ) + \sigma^\prime (\xi) F = 0   & \\[5pt]
\del_t v = \del_{xx} v + \del_x  \big (\overline{\sigma(u)} \big )    & \\[5pt]
S =  \overline{\sigma(u)}  + v_x   = \int \sigma(\xi) dF_{t,x} (\xi)  + v_x    & \\
\end{cases}
\quad (x,t) \in Q_T \, \; \; \xi \in \R
\end{equation}
subject to boundary and initial data
\begin{align}
S(0,t) &= S(1,t) = 0  \, , \quad   0 < t < T
\label{effbc1}
\\
v(0,x) &= v_0 (x) = \mbox{s-lim} v_0^\eps (x) \, , \quad F(0,x,\xi) = \mbox{wk-$\ast$ lim} \Big (  \charf_{ u_0^\eps (x) < \xi } \Big )
\label{effdata1}
\end{align}
where the functions $v(t,x)$ and $F(t,x,\xi)$ are defined in \eqref{wklim} and \eqref{defF}.

To see the first equation in \eqref{effs1}, start from \eqref{epsfamily} with $\beta(\cdot) \in C_c^1 (\R)$ and use the
representation formula \eqref{repform1} and $\del_\xi \charf_{u^\eps < \xi} = \delta (\xi - u^\eps)$ to write
$$
\begin{aligned}
\del_t \Big ( - \int \beta^\prime(\xi) \charf_{u^\eps < \xi} d\xi \Big ) &= \int \big (S^\eps - \sigma(\xi) \big )   \beta^\prime(\xi) \delta(\xi - u^\eps) d\xi
\\
&= \int \beta^\prime(\xi)  \Big [ \del_\xi \Big ( \big(S^\eps - \sigma(\xi) \big ) \charf_{u^\eps < \xi} + \sigma^\prime(\xi) \charf_{u^\eps < \xi} \Big ] d\xi
\end{aligned}
$$
From here we conclude that in the sense of distributions
\begin{equation}
\del_t \charf_{u^\eps < \xi}  + \del_\xi \Big (  \big (S^\eps - \sigma(\xi) \big ) \charf_{u^\eps < \xi} \Big ) + \sigma^\prime(\xi) \charf_{u^\eps < \xi} = 0
\end{equation}
Using \eqref{defF}, \eqref{stressb} and \eqref{stressconv} we pass to the limit $\eps \to 0$ and obtain \eqref{effs1}$_1$. The remaining
equations easily follow from \eqref{wklim}.

The equation \eqref{effs1}$_1$ may be written in the form
\begin{equation}
\del_t F + ( S-\sigma(\xi) ) \del_\xi F = 0
\end{equation}
and this writing is rigorous when $\sigma(\xi)$ is $C^1$ and $\del_\xi F$ is a measure. Let $f(u)$ be a $C^1$ function.
Using \eqref{moments}, we compute
\begin{align}
\del_t \overline{f(u)} &= \del_t \Big ( f(0) - \int_{-\infty}^0 f^\prime (\xi) F d\xi - \int_0^\infty f^\prime (\xi) (F-1) d\xi \Big )
\nonumber
\\
&= - \int_{-\infty}^0 f^\prime (\xi)  \del_t F d\xi - \int_0^\infty f^\prime (\xi) \del_t (F-1) d\xi
\nonumber
\\
&= \int_{-\infty}^0 (S - \sigma) f^\prime \del_\xi F d\xi + \int_0^\infty (S-\sigma) f^\prime \del_\xi F d\xi
\nonumber
\\
&= \int_{(-\infty, 0]} (S - \sigma) f^\prime  dF_{t,x} (\xi) + \int_{(0,\infty)} (S-\sigma) f^\prime dF_{t,x} (\xi)
\nonumber
\\[5pt]
&= \overline{(S-\sigma) f^\prime}
\label{formula3}
\end{align}

We apply this formula to the test function $f(u) = 0$ and obtain
$$
\del_t u = \overline{(S-\sigma) } = v_x
$$
which provides a derivation of \eqref{vemodel}$_1$ for solutions of the system \eqref{effs1}. One may also use \eqref{formula3}
to derive an alternate form of the effective system. Namely, using \eqref{formula3} for $f(u) = \sigma(u)$, we compute
$$
\del_t S = v_{xt} + \del_t  \overline{\sigma(u)} = S_{xx} + \overline{(S-\sigma) \sigma^\prime}
$$

This leads to the effective system
\begin{equation}\label{effs2}
\begin{cases}
\del_t F + \del_\xi \Big ( \big ( S - \sigma(\xi) \big ) F \Big ) + \sigma^\prime (\xi) F = 0   & \\[5pt]
\del_t S = \del_{xx} S + \int \sigma^\prime(\xi) (S - \sigma(\xi) ) dF_{t,x} (\xi)    & \\
\end{cases}
\; , \quad (x,t) \in Q_T \, \; \; \xi \in \R \, ,
\end{equation}
subject to boundary and initial data
\begin{align}
S(0,t) &= S(1,t) = 0  \, , \quad   0 < t < T \, , 
\label{effbc2}
\\
S(0,x) &= S_0 (x) = \mbox{wk lim}\big ( \sigma(u^\eps_0(x)) + \del_x v_0^\eps \big ) \, , \quad F(0,x,\xi) = \mbox{wk-$\ast$ lim} 
\Big (  \charf_{ u_0^\eps (x) < \xi } \Big ) \, .
\label{effdata2}
\end{align}
An alternative way for deriving \eqref{effs2}$_2$ is to start from \eqref{eqstress}  and pass to the limit $\eps \to 0$ using
\eqref{stressconv} and Young measures to represent the reaction term.

We summarize:

\begin{theorem}\label{mainthm}
Let $(u^\eps, v^\eps)$ and $S^\eps$ be a family of solutions of \eqref{vemodel}-\eqref{veid} induced by data satisfying \eqref{ub2data}
with $\sigma(u)$  nonmonotone  satisfying \eqref{boundW}, \eqref{hypgr}. Then, a subsequence defines $u(t,x)$, $v(t,x)$, $S(t,x)$ via \eqref{wklim}
and $F(t,x,\xi)$ via \eqref{defF} or via \eqref{defdistr}. The function $(F, v)$ satisfies the system \eqref{effs1}-\eqref{effdata1} or alternatively
$(F,S)$ is determined via \eqref{effs2}-\eqref{effdata2}.
\end{theorem}

As already seen,  given $(F,v)$ then $u$ is determined by the moment
\begin{equation}\label{formu}
u = - \int_{-\infty}^0 F(t,x,\xi) d\xi - \int_0^\infty (F(t,x,\xi) -1) d\xi
\end{equation}
and satisfies $u_t = v_x$. Alternatively, given $(F, S)$ determined via \eqref{effs2}-\eqref{effdata2}, 
the function $u$ is again determined by \eqref{formu}. We wish to define now $v(t,x)$
by the relation
$$
v_t = S_x \, ,  \quad v_x = S - \overline{\sigma(u)}
$$
Due to  \eqref{effs2}$_2$ and the formula $\del_t  \overline{\sigma(u)} = \overline{(S-\sigma)\sigma^\prime}$ (following from \eqref{formula3}),
these equations define $v$ and satisfy the desired equations $v_t = S_x$ and $u_t = v_x$.

\bigskip

{\bf Acknowledgement} I would like to thank {\sc Miguel Urbano} for helpful discussions
 and {\sc Denis Serre} for pointing out the references \cite{Serre91, Hillairet07} on homogenization.


\bigskip

{\bf Data availability} 
The data supporting the findings of this study are available within the paper.

{\bf Competing Interests}
The author has no relevant competing interests to declare.

{\bf Funding} 
Research supported by KAUST baseline funds, No BAS/1/1652-01-01



\begin{thebibliography}{10}

\bibitem{AB82}
G.~Andrews and J.~M. Ball.
\newblock Asymptotic behaviour and changes of phase in one-dimensional
  nonlinear viscoelasticity.
\newblock {\em J. Differential Equations}, 44(2):306--341, 1982.

\bibitem{b-Feireisl04}
E.~Feireisl.
\newblock {\em Dynamics of viscous compressible fluids}, volume~26 of {\em
  Oxford Lecture Series in Mathematics and its Applications}.
\newblock Oxford University Press, Oxford, 2004.

\bibitem{b-Folland99}
G.~B. Folland.
\newblock {\em Real analysis}.
\newblock Pure and Applied Mathematics (New York). John Wiley \& Sons, Inc.,
  New York, second edition, 1999.
\newblock Modern techniques and their applications, A Wiley-Interscience
  Publication.

\bibitem{Hillairet07}
M.~Hillairet.
\newblock Propagation of density-oscillations in solutions to the barotropic
  compressible {N}avier-{S}tokes system.
\newblock {\em J. Math. Fluid Mech.}, 9(3):343--376, 2007.

\bibitem{Hoff86}
D.~Hoff.
\newblock Construction of solutions for compressible, isentropic
  {N}avier-{S}tokes equations in one space dimension with nonsmooth initial
  data.
\newblock {\em Proc. Roy. Soc. Edinburgh Sect. A}, 103(3-4):301--315, 1986.

\bibitem{KLST23}
K.~Koumatos, C.~Lattanzio, S.~Spirito, and A.~E. Tzavaras.
\newblock Existence and uniqueness for a viscoelastic {K}elvin-{V}oigt model
  with nonconvex stored energy.
\newblock {\em J. Hyperbolic Differ. Equ.}, 20(2):433--474, 2023.

\bibitem{b-LSU68}
O.~A. Lady\v{z}enskaja, V.~A. Solonnikov, and N.~N. Ural'ceva.
\newblock {\em Linear and quasilinear equations of parabolic type}, volume Vol.
  23 of {\em Translations of Mathematical Monographs}.
\newblock American Mathematical Society, Providence, RI, 1968.
\newblock Translated from the Russian by S. Smith.

\bibitem{b-Lions98}
P.-L. Lions.
\newblock {\em Mathematical topics in fluid mechanics. {V}ol. 2}, volume~10 of
  {\em Oxford Lecture Series in Mathematics and its Applications}.
\newblock The Clarendon Press, Oxford University Press, New York, 1998.
\newblock Compressible models, Oxford Science Publications.

\bibitem{LPT94}
P.-L. Lions, B.~Perthame, and E.~Tadmor.
\newblock Kinetic formulation of the isentropic gas dynamics and {$p$}-systems.
\newblock {\em Comm. Math. Phys.}, 163(2):415--431, 1994.

\bibitem{LPS96}
P.-L. Lions, B.~Perthame, and P.~E. Souganidis.
\newblock Existence and stability of entropy solutions for the hyperbolic
  systems of isentropic gas dynamics in {E}ulerian and {L}agrangian
  coordinates.
\newblock {\em Comm. Pure Appl. Math.}, 49(6):599--638, 1996.

\bibitem{b-Perthame02}
B.~Perthame.
\newblock {\em Kinetic formulation of conservation laws}, volume~21 of {\em
  Oxford Lecture Series in Mathematics and its Applications}.
\newblock Oxford University Press, Oxford, 2002.

\bibitem{PT00}
B.~Perthame and A.~E. Tzavaras.
\newblock Kinetic formulation for systems of two conservation laws and
  elastodynamics.
\newblock {\em Arch. Ration. Mech. Anal.}, 155(1):1--48, 2000.

\bibitem{Serre91}
D.~Serre.
\newblock Variations de grande amplitude pour la densit\'{e} d'un fluide
  visqueux compressible.
\newblock {\em Phys. D}, 48(1):113--128, 1991.

\bibitem{Tzavaras23}
A.~E. Tzavaras.
\newblock Sustained oscillations in hyperbolic-parabolic systems.
\newblock {\em (preprint)}, 2023.

\end{thebibliography}

%
\end{document}